\documentclass[12pt]{amsart}
\usepackage{amsmath, amssymb,latexsym,amsthm, amscd}

\theoremstyle{definition}
\theoremstyle{plain}
\newtheorem{prop}[subsection]{Proposition}
\newtheorem{thm}[subsection]{Theorem}
\newtheorem{lem}[subsection]{Lemma}
\newtheorem{cor}[subsection]{Corollary}

\newcommand{\mbf}{\mathbf}
\newcommand{\mrm}[1]{\operatorname{#1}}

\title[Lie algebras and product Quivers]
{Generalized Kac-Moody Lie Algebras 
 And Product  Quivers}
\author{Yiqiang Li}
\address{Department of Mathematics\\ Kansas State University\\ 
Manhattan, KS  66506}
\email{yqli@math.ksu.edu}
\author{Zongzhu Lin}
\address{Department of Mathematics\\ Kansas State University\\ 
Manhattan, KS  66506}
\email{zlin@math.ksu.edu}
\thanks{The research was partially supported by NSF grant DMS-0200673}

\begin{document}

\begin{abstract}
We construct the entire generalized Kac-Moody Lie algebra as 
a quotient of the positive part of another generalized Kac-Moody Lie
algebra. 
The positive part of a generalized Kac-Moody Lie algebra 
can be constructed from representations of quivers using Ringel's
Hall algebra construction. 
Thus we give a direct realization of the entire generalized Kac-Moody
Lie algebra. 
This idea arises from the affine Lie algebra construction and
evaluation maps. 
In \cite{li-lin-C}, we give a
quantum version of this construction after analyzing Nakajima's quiver
variety construction of integral highest weight representations of the
quantized enveloping algebras in terms of the irreducible components of quiver varieties.

\end{abstract}

\maketitle

\section{Introduction}
Let $ \mrm{C}=(a_{ij})_{1\leq i, j\leq l}$ be a generalized Cartan matrix, and $
\mathfrak{g}(\mrm{C})$ be the Kac-Moody Lie algebra defined by generators $ e_i,
f_i$ and $ h_i$ ($i=1, \dots, l$) subject to the usual relations.
Then $ \mathfrak{g}(\mrm{C}) $ has a triangular decomposition 
\[ \mathfrak{g}(\mrm{C})=\mathfrak{n}^+\oplus \mathfrak{h}\oplus \mathfrak{n}^+.\]
We will use $ \mathfrak{n}^{+}(\mrm{C})$ etc, to indicate that the Lie algebras
are associated to the matrix $ \mrm{C}$. The main result of this paper is to prove
that for any $ \mrm{C}$ there exists another generalized Cartan matrix $ \mrm{C}'$ such
that $\mathfrak{g}=\mathfrak{n}^{+}(\mrm{C}')/I$ where $I$ has a specific set of
generators.

When $\mathfrak{g}(\mrm{C})$ is finite dimensional, it is clear from the
affinization construction of affine Kac-Moody Lie algebras that one can
simply take $ \mrm{C}'$ to be the un-twisted affine Cartan matrix. 
However, this
construction coincides with the construction in this paper only when
$ \mrm{C}=[2]$ (i.e., $ \mathfrak{g}(\mrm{C})=\mathfrak{sl}_2$). 

The motivation of this work was the construction of Kac-Moody Lie algebras
associated to symmetrizable generalized Cartan matrices by Peng and Xiao
\cite{peng-xiao} from the root category of the category of the representations of the quivers. 
It is known from the work of  Ringel that the Lie algebra 
$\mathfrak{n}^+(\mrm{C})$ can be constructed from the representations over finite
fields  of (valued) quiver $Q(\mrm{C})$
associated to $ \mrm{C}$ in terms of indecomposable representations. 
Later on,  in \cite{peng-xiao}, Peng and Xiao
went outside of the abelian category of representations of the quiver $ Q(\mrm{C})$, by
working in the root category of the abelian category (of representations of
$ Q(\mrm{C})$). Note that a root category is a quotient category of a derived
category, which is only a triangulated category with the translation functor
being a periodic functor of period 2. 
Thus, one can roughly see that two subalgebras 
$\mathfrak{n}^+(\mrm{C})$ and $ \mathfrak{n}^-(\mrm{C})$ arise from the abelian category
and its image under the translation functor.

After Peng and Xiao's work, a lot of effort have been made to construct the
entire quantum enveloping algebra $ U_q(\mrm{C})$ either from the root category or
any derived categories. But one crucial difficulty encounters:  the
obvious generalization of the Hall multiplication (the Lie bracket in
$\mathfrak{n}^+(\mrm{C})$ was realized in terms of Hall multiplication) fails to be
associative. The readers are referred to the work of Kapranov \cite{Kapranov2} and
To\"en \cite{Toen} for more discussions.

In this paper, instead of considering the representations of the
quiver $Q(\mrm{C})$, we consider a larger quiver $Q(\mrm{C})^{\pm}$ 
containing two copies of $Q(\mrm{C})$, 
each of its representations will give rise to a Lie subalgebra
$\mathfrak{n}^+(\mrm{C})$. 
Thus one will realize both $ \mathfrak{n}^+(\mrm{C})$ and  $
\mathfrak{n}^+(\mrm{C})$ in terms of the representations of the new quiver $
Q(\mrm{C})^{\pm}$. The motivation of using a larger quiver comes from Nakajima's
quiver variety construction in \cite{nakajima}, which provides basis to any
irreducible integrable highest weight modules of any (quantum) Kac-Moody Lie
algebra. In Nakajima's construction, the quiver $Q(\mrm{C})$ is enlarged to
duplicate the vertex set, but not the arrows set, linking the corresponding
duplicated vertices by two arrows of opposite direction. Of course Nakajima
uses the duplicated arrows to realize the fundamental weights. In our
construction, we work in the quiver of Nakajima's construction with arrows
of $ Q(\mrm{C})$ also duplicated to get $ Q(\mrm{C})^{\pm}$.
This construction for Lie
algebras also work for generalized Kac-Moody Lie algebras as well. 

The paper is organized as follows. 
In section 2, we recall some definitions.
In section 3, 
we give general construction of Lie algebras of the main result: 
For any generalized Kac-Moody Lie algebra $ \mathfrak{g}(\mrm{C})$ 
corresponding to the Borcherds-Cartan matrix $\mrm{C}$, 
we have a surjective Lie algebra homomorphism 
\[
\mathfrak{n}^+(\mrm{C}^{\pm})\rightarrow \mathfrak{g}(\mrm{C}),
\]
sending Chevalley generators to Chevalley generators. 
In section 4, 
we assume that $ \mrm{C}$ is symmetric and give a realization of
$\mathfrak{g}(\mrm{C})$ in terms of quiver representations of
$Q(\mrm{C})^{\pm}$. 
This is different from the construction of Peng-Xiao \cite{peng-xiao}.

{\bf Acknowledgements.} The authors thank Jie Xiao and Bangming Deng
for stimulating discussions on the questions of realizing the whole
Lie algebras.  In particular, Bangming Deng pointed out that there is
no need to introduce a new variable for the semisimple part.

\section{Generalized Kac-Moody Lie Algebras}
\subsection{}
Let $I$ be a countable index set.
Let $\mrm{\mrm{C}}=(a_{ij})_{i,j\in I}$ be 
a symmetrizable Borcherds-Cartan matrix, i.e., 
$\mrm{C}$ is an integer matrix  satisfying 
\begin{itemize}
\item[(BC1)] $a_{ii}=2$ or $a_{ii}\leq 0\quad$ for all $i\in I$,
\item[(BC2)] $a_{ij}\leq 0\quad $ if $ i\neq j$, and
\item[(BC3)] $a_{ij}=0$ implies $a_{ji}=0$
\end{itemize}
such that
there are nonnegative integers $\epsilon_i$ with no proper common divisor such that 
$\epsilon_i a_{ij}=\epsilon_j a_{j,i}$ for all $i, j \in I$.

We associate to  each Borcherds-Cartan matrix $\mrm{C}$
a generalized Kac-Moody Lie algebra  $ \mathfrak{g}(\mrm{C})$ 
with generators $e_i, h_i$ and $ f_i$ for all  $i \in I$ satisfying the following
relations.
\begin{itemize}
\item[(R1)] $[h_i,h_j]=0$.
\item[(R2)] $[h_i,e_j]=a_{ij}e_j$ and $[h_i,f_j]=-a_{ij}f_j$.
\item[(R3)] $[e_i,f_j]=\delta_{ij}h_i$.
\item[(R4)] For any $ i $ with $a_{ii}>0$ and all $ j\in I$
 \begin{itemize} 
 \item[(1)] $(\mrm{ad}(e_{i}))^{-2a_{ij}/a_{ii}+1}e_j=0,\quad $ and
 \item[(2)] $ (\mrm{ad}(f_{i}))^{-2a_{ij}/a_{ii}+1}f_j=0.$
 \end{itemize}
\item[(R5)] If $a_{ij}=0$, then $[e_i, e_j]=[f_i,f_j]=0$.
\end{itemize}
The Lie algebra $\mathfrak{g}(\mrm{C})$ has a triangular decomposition
\[ \mathfrak{g}(\mrm{C})=\mathfrak{n}^-\oplus \mathfrak{h}\oplus \mathfrak{n}^+\]
where $ \mathfrak{n}^-$ and $\mathfrak{n}^+$ are 
the  Lie subalgebras generated by $ f_i$'s  and $e_i$'s, respectively 
($i \in I$) 
and $\mathfrak{h}$ is the abelian subalgebra generated by $ h_i$'s.
When we are working with more than one Borcherds-Cartan matrices, we will
write $ \mathfrak{n}^+(\mrm{C})$, $\mathfrak n^-(\mrm{C})$ and $\mathfrak h(\mrm{C})$ 
to indicate that they are defined with respect to the matrix $\mrm{C}$. 
The relations in (R4) and (R5) are called Serre relations.

\subsection{} 
Define a Borcherds-Cartan matrix 
$\mrm{C}^{\pm}=(A_{(\alpha,i),(\beta,j})$, 
where $(\alpha,i),(\beta,j)\in \{+,-\}\times I$ as follows.
(For convenience, we write  $A^{\alpha,\beta}_{i,j}=A_{(\alpha,i),(\beta,j)}$.)
\[ A^{\alpha,\beta}_{i,j}
=\begin{cases} a_{ij} & \text{ if }
    \alpha=\beta,\\
-2\delta_{ij} & \text{ if }  \alpha\neq \beta.
\end{cases}\]
In  matrix form,  we  have 
$\mrm{C}^{\pm}=\begin{pmatrix}
    \mrm{C} & -2\mrm{Id}\\
  -2\mrm{Id} &   \mrm{C}\end{pmatrix}$,
where $\mrm{Id}$ stands for the identity $|I|\times |I|$ matrix. 
By definition, we know that 
$\mrm{C}^{\pm}$ is a symmetrizable
Borcherds-Cartan matrix.
In fact,
let $\epsilon^+_{i}=\epsilon_i$ and $\epsilon^-_i=\epsilon_i$ for all $i\in I$. Then the 
set $\{\epsilon^{\alpha}_i|(\alpha,i)\in \{+,-\}\times I\}$ is the minimal symmetrization 
of the Borcherds-Cartan matrix $\mrm{C}^{\pm}$.

The generalized Kac-Moody Lie algebra $\mathfrak{g}(\mrm{C}^{\pm})$ 
associated to $\mrm{C}^{\pm}$ 
has the generator set 
\[ \{ e^{+}_{i}, e^{-}_{i}, f^{+}_{i}, f^{-}_{i}, h^{+}_{i}, h^{-}_{i}\;|\;
i \in I\}\]
with the following generating relations:

\begin{itemize}
\item[(R1)] $[h^{\alpha}_i,h^{\beta}_j]=0\quad$ 
            for all $\alpha, \beta\in \{+,-\}$ and $i, j \in I$;
\item[(R2)] $[h^{\alpha}_i,e^{\alpha}_j]=a_{ij}e^{\alpha}_j\;\;$ and
            $\;\;[h^{\alpha}_i,f^{\alpha}_j]=-a_{ij}f^{\alpha}_j$, \\
            $[h^{\alpha}_i,e^{\beta}_j]=2\delta_{ij}e^{\beta}_j\;\;$ and
            $\;\;[h^{\alpha}_i,f^{\beta}_j]=-2\delta_{ij}f^{\beta}_j,\quad$ 
            for $\alpha\neq \beta$;
\item[(R3)] $[e^{\alpha}_i,f^{\alpha}_j]=\delta_{ij}h^{\alpha}_i$ and 
            $[e^{\alpha}_i,f^{-\alpha}_j]=0\quad $ $(i\neq j)$;
\item[(R4)] (1) $(\mrm{ad}(e^{\alpha}_{i}))^{-2a_{ij}/a_{ii}+1}e^{\alpha}_j
            =(\mrm{ad}(f^{\alpha}_{i}))^{-2a_{ij}/a_{ii}+1}f^{\alpha}_j=0$
            for any $ i $ with $a_{ii}=2$ and $ j\neq i$;\\
            (2) $(\mrm{ad}(e^{\alpha}_{i}))^{3}e^{-\alpha}_i
            =(\mrm{ad}(f^{\alpha}_{i}))^{3}f^{-\alpha}_i=0 \quad$ 
            for any $i$ with $a_{ii}=2$;\\            
            (3) $[e^\alpha_i, e^{-\alpha}_j]=[f^\alpha_i, f^{-\alpha}_j]=0\quad $ 
            for any  $ i\neq j$;
\item[(R5)] $[e^{\alpha}_i, e^{\alpha}_j]=[f^{\alpha}_i,f^{\alpha}_j]=0\quad $ 
            if $a_{ij}=0$.
\end{itemize}
Let $\mathfrak{n}^{-}(\mrm{C}^{\pm})$
(resp. $\mathfrak{n}^{+}(\mrm{C}^{\pm}))$ 
be the Lie subalgebra of $\mathfrak{g}(\mrm{C}^{\pm})$ 
generated by $f^{\alpha}_{i}$ (resp. $e^{\alpha}_{i}$),
where $i \in I$ and  $\alpha \in \{ +, -\}$.
Let $\mathfrak h(\mrm{C}^{\pm})$ be the Lie subalgebra 
of $\mathfrak g(\mrm{C}^{\pm})$ generated by $h_i^a$,
where $i\in I, \alpha \in \{+,-\}$. 
Then we have a standard  triangular decomposition
\[ 
\mathfrak{g}(\mrm{C}^{\pm})=\mathfrak{n}^{-}(\mrm{C}^{\pm})\oplus
\mathfrak{h}(\mrm{C}^{\pm})\oplus \mathfrak{n}^{+}(\mrm{C}^{\pm}).
\]
\subsection{} 
For each $i\in I\setminus I^{im}$, 
let $\mathfrak{g}(\mrm{C}^{\pm})_i$ be  the
Lie subalgebra of $ \mathfrak{g}(\mrm{C}^{\pm})$ generated by 
$f^{\alpha}_{i}, h^{\alpha}_i, e^{\alpha}_i$, isomorphic to the Kac-Moody
Lie algebra $ \widehat{\mathfrak{sl}}_2$. 

\section{$\mathfrak{g}(\mrm{C}) $ as a quotient of $\mathfrak{n}^{+}(\mrm{C}^{\pm})$}

\subsection{} 
\label{sl2}
Let $\mrm{C}=[2]$ be the Cartan matrix of $\mathfrak{sl}_2$.
Then $ \mrm{C}^{\pm}=\begin{pmatrix}2&-2\\-2&2\end{pmatrix}$. 
The generators of $ n^+(\mrm{C}^{\pm})$ are $ e^+, e^-$ with relations 
$\mrm{ad}(e^+)^3(e^-)=\mrm{ad}(e^-)^3(e^+)=0$.
 
Let $ \mathfrak{F}$ be the free Lie algebra of generators $x^+, x^-$
and let $ \mathfrak{J}=\ker{\phi}$ with 
$ \phi:\mathfrak{F}\rightarrow \mathfrak{sl}_2$ 
being the Lie algebra homomorphism defined by 
$ \phi(x^+)=e$ and $ \phi(x^-)=f$. 
Then $ \mathfrak{J}$ is generated by the following elements in $\mathfrak{F} $  
\begin{align}
\label{eq1}
[[x^+,x^-], x^+]-2x^+, \quad [[x^+,x^-], x^-]+2x^-.
\end{align}
If $ \psi: \mathfrak{F}\rightarrow \mathfrak n^+(\mrm{C}^{\pm})$ is the Lie
algebra homomorphism defined by $ \psi(x^+)=e^+$ and $ \psi(x^-)=e^-$,
then the ideal $  \mathfrak{I}=\ker{\psi}$ is generated by the two
elements
\begin{align}
\label{eq2}
\mrm{ad}(x^+)^3(x^-), \quad \mrm{ad}(x^-)^3(x^+).
\end{align}
By applying $ [-, x^+]$ and $[-,x^-]$ respectively to the two
generators of $ \mathfrak{J}$, we get the two generators of $
\mathfrak{I}$. Thus we have $\mathfrak{I}\subseteq \mathfrak{J}$. 

\begin{lem}
\label{sl2lem}
The map $e^+\mapsto e$, $e^-\mapsto f$ defines 
a surjective  Lie algebra homomorphism 
$ \pi: \mathfrak{n}^+(\mrm{C}^{\pm}) \to \mathfrak{sl}_2$ 
such that  $\phi=\pi \circ \psi$.  
The ideal  $\mathfrak{J}'=\ker(\pi) $ is generated by two elements  
\begin{align}
\label{eq3} 
[[e^+,e^-], e^+]-2e^+ \quad \text{and}\quad   
[[e^+,e^-], e^-]+2e^-.
\end{align} 
\end{lem}

\begin{proof}
In fact, let $\mathfrak k$ be the Lie ideal in $\mathfrak n^+(\mrm C^{\pm})$ generated by the elements in (\ref{eq3}).
Let 
\[
\pi': \mathfrak n^+(\mrm C^{\pm})/\mathfrak J' \to \mathfrak{sl}_2
\] 
be the induced Lie algebra isomorphism from $\pi$.
It is straight forward to check that $\mathfrak k \subseteq \mathfrak J'$.
Thus we have a surjective Lie algebra homomorphism
\[
\begin{CD}
\mathfrak n^+(\mrm C^{\pm})/\mathfrak k @>a>> \mathfrak n^+(\mrm C^{\pm})/\mathfrak J'.
\end{CD}
\]
On the other hand, by the universality of $\mathfrak{sl}_2$, there exists a Lie algebra homomorphism 
\[
\begin{CD}
\mathfrak{sl}_2 @>b>> \mathfrak n^+(\mrm C^{\pm})/\mathfrak k
\end{CD}
\]
sending $e \mapsto e^+$, $f \mapsto e^-$ and $h \mapsto [e^+,e^-]$.
Now by checking the generators, one sees that $b\circ \pi'$ is the inverse of $a$. Therefore $\mathfrak k=\mathfrak J'$.
\end{proof}

Note that one might think that $ \mathfrak{n}^{+}(\mrm{C}^{\pm})$ is
nilpotent as its standard name called as the nilpotent radical of
the Borel subalgebra of the affine Lie algebra. But $ \mathfrak{sl}_2$
is semisimple.  The above construction of $\pi$ is nothing special if
one think of $ \mathfrak{g}(\mrm{C}^{\pm})$ as the affinization of the
finite dimensional semisimple Lie algebra and then evaluate at 1, one
would also get a surjective homomorphism. However, the above
construction works for any generalized Kac-Moody Lie algebra. 

\subsection{} 
\label{general}
Let $ \mrm{C}$ be a Borcherds-Cartan matrix. Let $ \mathfrak{F} $ be
the free Lie algebra generated by $\{ x^{+}_i, x^{-}_i\;|\; i \in I\}$.
Let $ \phi: \mathfrak{F}\rightarrow \mathfrak{g}(\mrm{C})$ be the Lie algebra
homomorphism defined by $ \phi(x_i^+)=e_i $ and $ \phi(x^{-}_i)=f_i$.
The Lie ideal $ \mathfrak{J}=\ker(\phi)$ in $\mathfrak{F} $ is generated by 
\begin{equation} 
[[x^+_{i}, x_{i}^{-}],[x^+_{j}, x_{j}^{-}]], \quad
[[x^{-\alpha}_i, x^{\alpha}_i],x_j^{\alpha}]+a_{ij}x^{\alpha}_{j}, 
\quad ( i, j \in I, \alpha \in \{ +, -\}) 
\label{torus-rel} 
\end{equation}
and for $i\neq j \in I$
\begin{equation} \label{commut-rel}
[x^\alpha_i, x^{\alpha}_j] \quad \text{(if $ a_{ij}=0$)},\quad 
[x^\alpha_i, x^{-\alpha}_j] \quad (\alpha=+,-), 
\end{equation}
and for all $ i \in I$ with $ a_{ii}>0$ and $ j\neq i$
\begin{equation}\label{serre-rel}
\mrm{ad}(x^{\alpha}_{i})^{-a_{ij}+1}(x^{\alpha}_j) \quad (\alpha=+, -).
\end{equation} 
Similarly, 
let $ \psi: \mathfrak{F}\rightarrow \mathfrak{n}^{+}(\mrm{C}^\pm)$ 
be the Lie algebra homomorphism defined by 
$\psi(x_i^{+})=e_{i}^{+}$ and $ \psi(x_{i}^{-})=e^-_i$. 
The ideal $\mathfrak{I}=\ker{\psi}$ is generated by 
\begin{equation}
 \label{commut-rel+}
 [x^\alpha_i, x^{\alpha}_j] \quad \text{(if $ a_{ij}=0$)},\quad 
 [x^\alpha_i, x^{-\alpha}_j] \quad (i\neq j, \; \alpha=+,-), 
\end{equation}
and 
\begin{equation} \label{serre-rel+}  
 \mrm{ad}(x^{\alpha}_{i})^{-a_{ij}+1}(x^{\alpha}_j),\;
 \mrm{ad}(x^{\alpha}_{i})^{3}(x^{-\alpha}_i) \;\;\quad 
 ( a_{ii}>0, \; i\neq j,\; \text{ and }\alpha=+, -). 
\end{equation} 
Note that the relations in \eqref{commut-rel} and  \eqref{commut-rel+} are
the same. The only difference between \eqref{serre-rel} and \eqref{serre-rel+} is
the elements $\mrm{ad}(x^{\alpha}_{i})^{3}(x^{-\alpha}_i) $ 
for all $i$ with $ a_{ij}>0$. But 
\[ \mrm{ad}(x^{\alpha}_{i})^{3}(x^{-\alpha}_i)=[[[x^{-\alpha}_i,
x^{\alpha}_i],x_i^{\alpha}]+a_{ii}x^{\alpha}_{i}, x^{\alpha}_{i}] \in
\mathfrak{J}.\]
Thus we have $ \mathfrak{I}\subseteq \mathfrak{J}$. 

\begin{thm}
\label{mainthm}
There is a surjective 
Lie algebra homomorphism of 
$\pi: \mathfrak{n}^+(\mrm{C}^\pm)\rightarrow \mathfrak{g}(\mrm{C})$.
Moreover, the kernel $ \ker{(\pi)}$ is the ideal of $ \mathfrak{n}^{+}(\mrm{C}^\pm)$
generated by 
\begin{equation}
\label{eq9} 
[[e^+_{i}, e_{i}^{-}],[e^+_{j}, e_{j}^{-}]], \quad
[[e^{-\alpha}_i, e^{\alpha}_i],e_j^{\alpha}]+a_{ij}e^{\alpha}_{j}, 
\quad ( i, j \in I, \alpha \in \{ +, -\}).
\end{equation}
\end{thm}

The proof of Theorem ~\ref{mainthm} is the same as the proof of Lemma ~\ref{sl2lem}.
First of all, let $\mathfrak k$ be the Lie ideal of $\mathfrak n^+(\mrm C^{\pm})$ generated by elements in (\ref{eq9}). Then
by the universality of $\mathfrak g(\mrm C)$, we have a surjective Lie algebra homomorphism 
$\mathfrak g(\mrm C) \to \mathfrak n^+(\mrm C^{\pm})/\mathfrak k$. 
This is in a sense the inverse of the map 
$\mathfrak n^+(\mrm C^{\pm})/\mathfrak k \to \mathfrak n^+(\mrm C^{\pm})/\ker{(\pi)}$.
See Proof of Lemma ~\ref{sl2lem} for more details.

\section{Degenerate Ringel-Hall Algebras}

\subsection{} 
In this section, 
we take $ Q=(I, \Omega, s, t)$ as an arbitrary locally finite quiver,
i.e.,  finite vertex subset generate a finite subquiver. 
Let us fix a prime $p$. 
For any $r$, 
let $ \mathbb{F}_r =\mathbb{F}_{p^r}$ the finite field of $p^r$ elements,
which can be regarded as a subfield of 
$\bar{\mathbb{F}}$, the algebraic closure of  $\mathbb{F}_p$.

Let $ \mathcal{M}(Q, r)$ the abelian category of all finite dimensional
nilpotent representations of $ Q$ over $\mathbb{F}_r$. 
If $ M$ is an object of $ \mathcal{M}(Q,r)$, 
we will denote by $M(r')=M\otimes_{\mathbb{F}_r}\mathbb{F}_{r'}$ 
for any $r' =mr$, a multiple of $ r$.

Let $H(Q, r)$ be the free abelian group with basis being the isomorphism
classes of objects in $\mathcal{M}(Q,r)$. For each object $M$ in
$\mathcal{M}(Q,r)$, let $ d(M)$ the number of components in the
decomposition of $ M$ into direct sum of indecomposable representations. 
Then $ H(Q,r)$ is graded by $\mathbb{N}$ with the $d$-component $ H(Q,r)_d$
being the free abelian module generated by isomorphism classes $ [M]$ with $
d(M)=d$. 

Ringel defined a comultiplication 
$ \Delta: H(Q,r)\rightarrow H(Q,r)\otimes_{\mathbb{Z}} H(Q,r)$ 
by 
\[  
\Delta([M])=\sum_{([M_{1}], [M_{2}])} [M_1]\otimes [M_2],
\]
where $([M_1], [M_2])$ runs over all pairs of isomorphism classes such
that $[M_1 \oplus M_2]=[M]$.
\begin{lem}[Ringel \cite{Ringel:liealg}] $H(Q,r)$ is a strictly graded
 cocommutative   $\mathbb{Z}$-coalgebra, with co-unit defined by
 $\epsilon([M])=\delta_{[0], [M]} $.
\end{lem}
On the other hand, $H(Q, r)$ has a natural associative  $
\mathbb{Z}$-algebra structure with  Hall multiplication 
\[[M][N]=\sum_{[E]} h_{M, N}^{E} [E],\]
where $  h_{M, N}^{E} $ is the Hall number, i.e., the number of
subrepresentatins 
$K$ of $E$ such that $[K]=[N] $ and $ [E/K]=M$ and $ [0]$ is the identity element. 
With respect to the above algebra and coalgebra structure, $ H(Q,r)$
is not a bialgebra in general. But we have  
 
\begin{lem}[Ringel \cite{Ringel:liealg}]
With the natural induced algebra and coalgebra structure on 
$\bar{H}(Q,r)=H(Q, r)\otimes_{\mathbb{Z}}\mathbb{Z}/(p^r-1)$,
$\bar{H}(Q,r)$ is a bialgebra over the commutative ring $
R_r=\mathbb{Z}/(p^r-1)$
\end{lem}
\begin{proof}
This is a consequence of \cite[Prop.2, Prop.3]{Ringel:liealg}.
\end{proof}

Recall that the set of positive integers is 
a poset under the partial order $r|r'$. 
Any pair of positive integers has a common multiple.  
Note that if $r|r'$, then $ (p^r-1)|(p^{r'}-1)$. 
Thus there is a natural ring homomorphism 
$ R_{r'}\rightarrow R_r$. 
So we have an inverse system of finite rings
over the poset of all positive integers under the divisibility order. 
Let $ \bar{R}=\varprojlim_{r} R_r$, which contains $ \mathbb{Z}$ as subring. 
Note that $\bar{R}$ is a subring of  $\prod_{r>0} R_r$.
$ \mathbb{Z}\rightarrow \prod_{r>0}R_r$ is an injective ring
homomorphism. Then $ \prod_{r>0}\bar{H}(Q,r)$ is an associative
$\bar{R}$-algebra with the componentwise multiplication.

Although each $ \bar{H}(Q,r)$ is a bialgebra, the comultiplication does not
naturally extend to $\prod_{r>0}\bar{H}(Q,r)$.

\subsection{} 
For two elements 
$(u_r)_{r>0},  (v_r)_{r>0} \in \prod_{r>0}\bar{H}(Q,r)$, 
we define $ (u_r)\sim (v_r)$ if there exists $r_0>0$ such that 
$u_r=v_r$ for all $r$ with $r_0|r$. 
This is an equivalent relation and is compatible with 
the $\bar{R}$-algebra structure, i.e.,
if $ (u_r)\sim (u'_r)$, $ (v_r)\sim (v_r')$, and $ a=(a_r)\in \bar{R}$,
then 
\begin{eqnarray}
 (u_r)+(v_r)\sim (u'_r)+(v'_r), & (u_r)(v_r)\sim (u'_r)(v'_r) & \text{and} \quad a(u_r)\sim a(u'_r).
\end{eqnarray} 

Let $\bar{\Pi}=\prod_{r>0}\bar{H}(Q,r)/\sim$ be the set of equivalence classes
which has a natural $ \bar{R}$-algebra structure such that the quotient map
$ \prod_{r>0}\bar{H}(Q,r)\rightarrow \bar{\Pi}$ is a homomorphism of 
$\bar{R}$-algebra.

We will write the image of $ (u_r)\in \prod_{r>0}\bar{H}(Q,r)$ still by $ (u_r)$.

For each $r_0>0$ and  each object $M$ in $\mathcal{M}(Q,r_0)$, we define
$\{M\}=([M_r])\in \prod_{r>0}\bar{H}(Q,r)$ by 
\[ 
[M_r]=\begin{cases} [M(r)] & \text{ if $r_0 |r $ }\\
[0] & \text{ otherwise }.\end{cases}
\] 
It follows from the definition that we have $ \{M\}\neq \{M(r)\}$,
if $ r_0\neq r$,  in $\prod_{r>0} \bar{H}(Q,r)$. 
However, we have $ \{M\} \sim \{M(r)\}$. 
Thus $\{M\}$ is a well-defined element in $ \bar{\Pi}$,
independent of the actual field representation of $ M$.

For each $ M \in \mathcal{M}(Q,r_0)$, 
one can choose a finite field extension $\mathbb{F}_{p^r}$ 
such that all indecomposable components of $ M(r)$ are
absolutely indecomposable. 
Then 
\[
\Delta(\{M\})=\big (\sum_{[M_1\oplus M_2]=[M]} [M_{1}]\otimes
[M_{2}] \big )
=\sum_{[M_1(r)\oplus M_2(r)]=[M(r)]} \{M_1(r)\} \otimes \{M_2(r)\}
\]
is a well-defined element in $\bar{\Pi}$. 
(Note that the right hand side is a finite sum.)

Let $ H_1(Q)$ be the  $\bar{R}$-subalgebra of $ \bar{\Pi}$ generated
by all $\{M\}$ with  $ M $ in $ \mathcal{M}(Q,r)$,  for  all $ r>0$.

\begin{lem} 
$\Delta$ extends uniquely to a comultiplication making $ H_1(Q)$ a bialgebra. 
\end{lem}
  
\begin{proof} 
On each generators of $H_1(Q)$, $\Delta$ is defined componentwisely. 
on each component, $\Delta$ is a homomorphism of algebras. 
Thus $ \Delta$ extends uniquely to $ H_1(Q)$ and making $H_1(Q)$ a bialgebra.  
\end{proof}

\begin{cor}
The set $L_1(Q)$ of primitive elements of $ H_1(Q)$ is an $ \bar{R}$-Lie algebra.
\end{cor}

Note that an element $x$ is called primitive if 
$\Delta(x)=x\otimes 1 + 1 \otimes x$.
Thus, $\{M\}$ is a primitive element if $M$ is  absolutely indecomposable.

{\bf Question:} Is $L_1(Q)$ generated by all $\{M\}$'s with $M$
absolutely indecomposable?

Let $L(Q)$ be the Lie subalgebra of $L_1(Q)$ generated by 
$\{ S_i\}$ ($i\in I$), where $S_i$ is the simple nilpotent representation
associated to the vertex $i$.  
Let $\mrm C=(a_{ij})_{i,j\in I}$ be the (symmetric) Borcherds-Cartan matrix associated to 
$Q$. 
It is defined by $a_{ii}=2(1-l_i)$ and $-a_{ij}$ = the number of arrows
between $i$ and $j$, if $i\neq j$, where $l_i$ is the number of loops
at $i$.
Recall that $\mathfrak g(\mrm C)$ is the generalized Kac-Moody Lie
algebra associated to $\mrm C$. 
$\mathfrak n^+(\mrm C)$ is its positive part.
We have 

\begin{prop}
There is a unique Lie algebra isomorphism $\phi: \mathfrak n^+(\mrm C) \to L(Q)$
such that $\phi(e_i)=(S_i)$, for all $i\in I$.
\end{prop}

In fact, this is a consequence of the results in 
\cite{Ringel:hallqua, Green, KS}, when the quiver is of finite type,
with no loops and with possible loops, respectively.

\subsection{} 
Let $Q'=(I', \Omega', s', t')$ be the Kronecker quiver. 
Specifically, $I'$ is the set $\{+, -\}$ and  
$\Omega'$ is the set $\{\alpha,\beta\}$
such that $s'(\alpha)=+$, $t'(\alpha)=-$,
$s'(\beta)=+$ and $t'(\beta)=-$.

Given any quiver $Q=(I, \Omega, s, t)$, we define the quiver
$Q^{\pm}=(I', \Omega', s', t')$
as follows.

\begin{itemize}
 \item[(1)] $I''=I\times I' $.
 \item[(2)] $H''=I\times H' \cup H\times I'$ (disjoint union).  
 \item[(3)] $s''((a,b))=
             \left\{
              \begin{array}{ll}
               (a,s'(b)) & \mbox{if $(a,b)\in I\times H'$},\\
               (s(a),b) & \mbox{if $(a,b)\in H\times I'$}.
              \end{array}
             \right.$
 \item[(4)] $t''((a,b))=
             \left\{
              \begin{array}{ll}
                (a,t'(b)) & \mbox{if $(a,b)\in I\times H'$},\\
                (t(a), b) & \mbox{if $(a,b)\in H\times I'$}.
              \end{array}
              \right.$
\end{itemize}  
We call $Q^{\pm}$ the product quiver of $Q$ and $Q'$. 
Of course, this definition works for any two quivers.
One can check that the Borcherds-Cartan matrix associated to $Q^{\pm}$
is exactly $\mrm C^{\pm}$, 
where $\mrm C$ is the Borcherds-Cartan matrix associated to $Q$.

Let $L(Q^{\pm})$ be the Lie algebra obtained by applying the
construction in the previous subsection. 
Note that $L(Q^{\pm})$ is a Lie algebra generated by
$\{S_{(i,a)}\}$, for all $(i,a)\in I''$, where $S_{(i,a)}$ is  
the simple nilpotent representation of $Q^{\pm}$ associated to the
vertex $(i,a)$.
For convenience, we simply write $S_{(i,a)}$ for $\{S_{(i,a)}\}$ when
no confusion will be made.
Applying Proposition 4.7 to the quiver $Q^{\pm}$, we have 
\begin{cor}
The assignment $e_i^a \mapsto S_{(i,a)}$ $((i,a)\in I'')$ defines a Lie algebra
isomorphism $\phi: n^+(\mrm C^{\pm}) \to L(Q^{\pm})$.
\end{cor}

By applying Theorem 3.3, we have 
\begin{thm}
There is a Lie algebra isomorphism 
\[
\pi: L(Q^{\pm})/\mathfrak i \to \mathfrak g(\mrm C),
\]
such that $\pi(S_{(i,+)})=e_i$ and $\pi(S_{(i,-)})=f_i$,
where $\mathfrak i$ is the Lie ideal of $L(Q^{\pm})$ generated by 
\begin{equation} 
[[S_{(i, +)}, S_{(i,-)}],[S_{(+,j)}, S_{(j,-)}]], 
[[S_{(i, -a)}, S_{(i,a)}],S_{(j,a)}]+a_{ij}S_{(j,a)}, 
\quad ( i, j \in I, a \in \{ +, -\}).
\end{equation}
\end{thm}

\noindent
{\bf Remark.}
When $\mrm C$ is symmetrizable,
one can realize $\mathfrak g(\mrm C)$ by considering the category of 
representations of the valued quiver
$\overset{\rightarrow}{\Gamma^{\pm}}$,
defined in [LL].

\subsection{Relation with the affinization}

Let $\mathfrak g=\mathfrak{sl}_2=\mathbb C\{e,f,h\}$.
From Section 7.6 in \cite{Kac}, we know that the affinization of $\mathfrak g$
is $\widehat{\mathfrak{sl}_2}$ and the positive part of 
$\widehat{\mathfrak{sl}_2}$ is 
\[
\mathfrak n^+= 
t\mathbb C[t]\otimes \mathbb C\{f, h\}
+
\mathbb C[t]\otimes \mathbb C\{e\},
\]
with the Lie bracket defined by 
\[
[t^m\otimes x, t^n\otimes y]=t^{m+n}\otimes [x,y],
\]
where $t$ is an indeterminate, $x,y\in \mathfrak g$, 
and $m,n$ are suitable nonnegative integers.
In particular, we have 
\begin{enumerate}
\item[(a)] $t\otimes h=[1\otimes e, t\otimes f]$;
\item[(b)] $t^{n+1}\otimes e=\frac{1}{2}[t\otimes h, t^n\otimes e]$;
\item[(c)] $t^{n+1}\otimes f=\frac{1}{2}[t^n\otimes f, t\otimes h]$;
\item[(d)] $t^n\otimes h=[1 \otimes e, t^n\otimes f]$.  
\end{enumerate}
From \cite{Kac}, we also know that $\mathfrak n^+$ is generated by
$e_0=t\otimes f$ and $e_1=1\otimes e$.

From Section 3.1, we know that there is a unique surjective Lie algebra
homomorphism $\pi_1: \mathfrak n^+ \to \mathfrak{sl}_2$ such that
$\pi_1(e_1)=e$ and $\pi_1(e_0)=f$. 
Moreover $\mrm{ker} \pi_1$ is a Lie ideal of $\mathfrak n^+$ generated
by 
\[
[[e_1, e_0],e_1]-2e_1, \quad 
[[e_1,e_0], e_0]+2e_0.
\]
By direct computation, we have
\[
[[e_1, e_0],e_1]-2e_1=(t-1) 2e_1, \quad
[[e_1,e_0], e_0]+2e_0=-(t-1) 2e_0.
\]
This implies that $\mrm{ker} \pi_1$ is a Lie ideal of $\mathfrak n^+$
generated by $(t-1)$.
Therefore, $\pi_1$ is nothing but the evaluation map, taking $t$ to
$1$.

On the other hand,
when the quiver $Q$ has only one vertex and with no arrows, 
the quiver $Q^{\pm}$ is the Kronecker quiver.
The indecomposable representations of $Q^{\pm}$ over $\mathbb F_r$
is classified into
three classes: preprojective, regular and preinjective. 
The dimension vectors of the indecomposable preprojectives, 
(resp. regulars, preinjectives) are of the form
$(n-1, n)$ (resp. $(n,n)$, $(n, n-1)$).
We denote by $S_+$ and $S_-$ the simple representations of $Q^{\pm}$
over $\mathbb F_r$
corresponding to the vertices $+$ and $-$, respectively. 
We also denote by  $P_n$ (resp. $I_n$)  the preprojective indecomposable of
dimension vector $(n-1, n)$ (resp. $(n,n-1)$). 
Let $R_n$ be the sum of all regular indecomposable representations  of dimension
vector $(n,n)$.
Note that $S_-=P_1$ and $S_+=I_1$.

In $H(Q,r)$, given any $x=\sum_i x_i[M_i]$, we set 
$Z(x)=\sum_{M_i\; \mbox{is regular}} x_i[M_i]$, the regular part of $x$.
By \cite{Z}, we have in $H(Q,r)$ the following relations:
\begin{enumerate}
\item $R_1=S_+S_--S_-S_+$;
\item $I_{n+1}=\frac{1}{q+1}(R_1I_n - q I_n R_1)$;
\item $P_{n+1}=\frac{1}{q+1}(P_nR_1 - q R_1 P_n)$;
\item $Z(I_jP_i)=I_jP_i-q^{i+j-2}P_iI_j$;
\item $Z(I_jP_i)=Z(I_1P_{i+j-1})=Z(I_{i+j-1}P_1)$;
\end{enumerate}
where $q=p^r$. 
By specializing $q$ to $1$, we have in $\bar H(Q,r)$ the following
relations:
\begin{enumerate}
\item[(1')] $R_1=[S_+, S_-]$;
\item[(2')] $I_{n+1}=\frac{1}{2}[R_1, I_n]$;
\item[(3')] $P_{n+1}=\frac{1}{2}[P_n, R_1]$;
\item[(4')] $Z(I_jP_i)=[I_j, P_i]$; 
\item[(5')] $Z(I_jP_i)=Z(I_1P_{i+j-1})=Z(I_{i+j-1}P_1)$.
\end{enumerate}
Moreover, since $p-1$ divides $h^E_{I_j, P_i}-h^E_{P_i, I_j}$ whenever
$E$ is a decomposable representation of $Q^{\pm}$, we have that 
$Z(I_jP_i)=\sum_{[M]\in R_{i+j-1}}h^M_{I_j, P_i} [M]$ in $\bar H(Q,r)$.
Given any  $[M] \in R_{i+j-1}$, one can check that
$h^M_{I_1, P_{i+j-1}}=\frac{q^{i+j}-1}{q-1}=\sum_{l=0}^{i+j-2}q^l$.
Thus, by $\mrm{(5')}$ we have in $\bar H(Q, r)$
\begin{enumerate}
\item[(6')] $Z(I_jP_i)=(i+j-1)R_{i+j-1}$.
\end{enumerate}
Therefore, in $L(Q^{\pm})$, we have the following relations:
\begin{enumerate}
\item[(i)]    $\{R_1\}=[\{S_+\}, \{S_-\}]$;
\item[(ii)]   $\{I_{n+1}\}=\frac{1}{2} [\{R_1\}, \{I_n\}]$;
\item[(iii)]  $\{P_{n+1}\}=\frac{1}{2}[\{P_n\}, \{R_1\}]$;
\item[(iv)]   $n\{R_n\}=[\{I_1\}, \{P_n\}]$; 
\end{enumerate}
where, by abuse of notation, $\{R_n\}$ is defined by
$\{R_n\}_r=R_n$, for any $r$.
(Note that $(\mrm{iv})$ comes from $\mrm{(4')-(6')}$.)

Now the assignments $e_1=1\otimes e\mapsto \{S_+\}$ and 
$e_0=t\otimes f \mapsto \{S_-\}$ define a Lie algebra isomorphism 
$\phi: \mathfrak n^+ \to \mathbb C \otimes_{\mathbb Z} L(Q^{\pm})$, 
by Proposition 4.7.
Moreover, by comparing the relations $\mrm{(a)-(d)}$ with the relations
$\mrm{(i)-(iv)}$, we have
\begin{enumerate}
\item[(I)]   $\phi: t^{n-1}\otimes e \mapsto \{I_n\}$;
\item[(II)]  $\phi: t^n\otimes f \mapsto \{P_n\}$;
\item[(III)] $\phi: t^n\otimes h \mapsto  n \{R_n\}$.
\end{enumerate}
From $\mrm{(I)-(III)}$, one has that
in $L(Q^{\pm})/\mathfrak i$,  
$\{S_-\}= \{P_n\}$, $\{ S_+\}=\{I_n\}$
and $\{R_1\}=m\{R_m\}$, for any $m, n$.

\end{document}